\documentclass[12pt]{article}
\usepackage{amsfonts,amssymb,amsmath}
\usepackage[T2A]{fontenc}  
\usepackage[cp1251]{inputenc}
\usepackage[english,russian]{babel}
\usepackage{geometry}
\newtheorem{theorem}{Theorem}
\newtheorem{lemma}[theorem]{Lemma}

\newtheorem{rem}[theorem]{Remark}

\def\ds{\displaystyle}
\def\f{\frac}
\def\div{{\rm div}}
\def\C{\mathbb{C}}
\def\R{\mathbb{R}}

\def\K{\mathcal{K}}

\def\per{{\rm per}}
\def\nab{\nabla}

\def\a{\alpha}
\def\b{\beta}

\def\ve{\varepsilon}
\def\e{\varepsilon}

\def\C0{C^\infty_0(\R^d)}

\def\ld{L^2(\R^d)}
\def\rd{\R^d}

\def\ild{\int_{\R^d}}
\def\ilb{\int_{Y}}
\def\beq{\begin{equation}}
\def\eeq{\end{equation}}

\begin{document}
	
	%
	%
	%
	%
	%
	%
	%
	%
	%
	
\title{Homogenization estimates\\
for high order elliptic
operators 
}


\author{ S.\,E. Pastukhova}

\date{}
\maketitle


\begin{footnotesize}

%
%
%
\noindent
In the whole space $\rd$, $d\ge 2$, we study homogenization of a divergence form elliptic operator
$A_\e$ of order $2m\ge 4$ with measurable $\varepsilon$-periodic coefficients, where
$\varepsilon$ is a small parameter. 
For the resolvent $(A_\e+1)^{-1}$, we construct an approximation
with the remainder term of order $\varepsilon^2$
in the operator
$(L^2{\to}H^m)$-norm, using the resolvent of the homogenized operator, solutions of several auxiliary periodic problems on the unit cube,  and smoothing operators. The homogenized operator here 
differs from the one commonly employed in homogenization.
\end{footnotesize}
	
	\section{Introduction}\label{sec:sec1}
	\setcounter{theorem}{0} 
\setcounter{equation}{0}
	
	The present paper is devoted to the homogenization of high order elliptic operators with periodic coefficients. We mean 
	operators of an arbitrary even order $2m\ge 4$.
	The first qualitative results in this topic were obtained long ago in 70s (see, e.g., \cite{B} and \cite{ZKOK}). In this paper, we are interested in estimates for the homogenization error, which are of operator type,  and continue the recent studies of \cite{P16}--\cite{P20s}, where the approach proposed in \cite{{Zh1}} and \cite{{ZhP05}} was applied in different situations concerning high order elliptic operators. Here, we construct  approximations for the resolvents of high order operators with the remainder term 
	$O(\varepsilon^2)$ as $\e\to 0$
	in the energy operator norm, i.e., in the operator
$(L^2{\to}H^m)$-norm.
	
{\bf 1.1.} In the whole space $\rd$, $d{\ge} 2$, we consider the following equation of an even order $2m{\ge}4$: 
\begin{equation}\label{1}
  \begin{array}{cc}
  u^\e\in  H^m(\rd),  \quad     (A_\e+1) u^\e
  =     f,  \quad f\in L^2(\rd),& \\
A_\e =
(-1)^m\sum\limits_{|\alpha|=|\beta|=m}D^{\alpha}(a^\e_{\alpha\beta}(x)D^{\beta}), & \\
  \end{array}
\end{equation}
with rapidly 
oscillating  $\varepsilon$-periodic coefficients  
$
a^\e_{\alpha\beta}(x)=a_{\alpha\beta}(y)|_{y=\e^{-1}x},
$
for small $\ve{\in}( 0,1)$.
Here, $D^{\alpha}$ denotes the multiderivative
\[
D^{\alpha}=D_1^{\alpha_1}\ldots D_d^{\alpha_d},\quad D_i=\f{\partial}{\partial x_i}, \quad i=1,\ldots,d,
\] where 
$
\alpha{=}(\alpha_1,\ldots,\alpha_d)$  is the multiindex of length
$|\alpha|{=}\alpha_1+\ldots+\alpha_d$ with $ \alpha_j{\in}
\mathbb{Z}_{\ge 0}$; the coefficients
$a_{\alpha\beta}(y)$ are measurable periodic real-valued functions, and
the unit cube $Y{=}[-1/2,1/2)^d$ is the periodicity cell. 
 
 We assume that the following  boundedness and ellipticity conditions hold:
\begin{equation}\label{2}\ds{
\|a_{\alpha\beta}\|_{ L^\infty(Y)}\le \lambda_1,}
\atop\ds{
 \int_{\rd}\sum\limits_{|\alpha|=|\beta|=m}a_{\alpha\beta}(x)D^{\beta}\varphi
 D^{\alpha}\varphi\,dx\ge  \lambda_0
\int_{\rd} \sum\limits_{|\alpha|=m}| D^{\alpha}\varphi|^2\,dx\quad \forall \varphi\in C_0^\infty(\rd)
}
\end{equation}
for some positive constants $ \lambda_0$, $\lambda_1$ and all multiindices $\alpha$ and $\beta$ of length $m$. No symmetry conditions are imposed on the coefficients $a_{\alpha\beta}$.

In (\ref{1}), we use the Sobolev space  $H^m{=}H^m(\rd)$, 
equipped with the norm 
 \[
\|u\|^2_{H^m}=\ild
\sum\limits_{|\alpha|\le m}|D^\alpha u|^2\,dx.
\]
As known, the set  $C_0^\infty(\rd)$ of smooth compactly supported functions is dense in $H^m(\rd)$  and the norm can be equivalently introduced by a simpler way 
\[
\|u\|^2_{H^m}=\ild
\sum\limits_{|\alpha|= m}|D^\alpha u|^2\,dx+
\ild
|u|^2\,dx,
\]
  
 $G$-convergence and homogenization issues for differential operators $A_\e$ in (\ref{1}) have been studied since the 70s.
Even more general operators were considered from this point of view in \cite{ZKOK}, namely, 
divergence form operators
\begin{equation}\label{3}
A_\e =\sum\limits_{|\alpha|\le m,|\beta|\le m}(-1)^{|\alpha|}D^{\alpha}(a^\e_{\alpha\beta}(x)D^{\beta}),
\end{equation}
 with lower order terms. Moreover, the periodicity of coefficients was not necessarily required for (\ref{3}), instead of which the so called $N$-condition was assumed in \cite{ZKOK}.
 The well known result on homogenization of the operator  $A_\e$ in (\ref{1}) means the closeness, for sufficiently small  $\varepsilon$, in the sense of the strong operator topology between the resolvent $(A_\e+1)^{-1}$, regarded as an operator in $ L^2(\rd)$, and the resolvent $(\hat A+1)^{-1}f$ of the homogenized operator $\hat A$.
 More exactly,
   $(A_\e+1)^{-1}f$ converge to
$(\hat A+1)^{-1}f$ in the $ L^2(\rd)$-norm as $\e\to 0$ for any    $f\in L^2(\rd)$.  The homogenized operator $\hat A$ is of the same class (\ref{2}) as the original operator $A_\e$, but much simpler:
  \begin{equation}\label{4}
 \hat A=(-1)^m
\sum\limits_{|\alpha|=|\beta|= m}D^{\alpha}\hat{a}_{\alpha\beta}D^{\beta},
\end{equation} 
where the coefficients  $\hat{a}_{\alpha\beta}$ are constant and can be 
defined with the help of auxiliary problems on the periodicity cell $Y$ (see Subsection 3.1).    The homogenized problem for (\ref{1}) is written as
 \begin{equation}\label{5}
   u\in  H^m(\rd),  \quad     \hat A u+u=
     f,  \quad f\in L^2(\rd),
\end{equation}
and the above result about the strong resolvent convergence means that the solutions to the problems
 (\ref{1}) and (\ref{5}) are connected by the limit relation 
\begin{equation}\label{5a}
\lim_{\e\to 0}\|u^\e-u\|_{L^2(\rd)}=0
\end{equation}
 for any right-hand side function $f\in L^2(\rd)$.

{\bf 1.2.}  In \cite{V}, \cite{P16} and \cite{AA16}, 
 the stronger operator convergence of $A_\e$ to $\hat A$
was established, namely, the uniform resolvent convergence in the operator $ L^2(\rd)$-norm; moreover, the convergence rate estimate with respect to the parameter $\varepsilon$ was deduced:
 \begin{equation}\label{6}
\|(A_\e+1)^{-1}-(\hat A+1)^{-1}\|_{L^2(\rd)\to L^2(\rd)}\le C\e,
 \end{equation}
where the constant $C$ depends only on the spatial dimension $d$ and the constants  $\lambda_0$ and $\lambda_1$  in (\ref{2}). The main result in 
\cite{P16} и \cite{AA16} is even more exact and stronger; it concerns approximation of the resolvent
  $(A_\varepsilon+1)^{-1}$
 in the operator $(L^2(\rd)\to H^m(\rd))$-norm by using the sum  $(\hat A+1)^{-1}+\e^m\mathcal{K}_\e$ of the resolvent of the homogenized operator $\hat A$ and the correcting operator. Furthermore,
\begin{equation}\label{7}
\|(A_\varepsilon+1)^{-1}-(\hat A+1)^{-1}-\varepsilon^m \mathcal{K}_\varepsilon\|_{L^2 (\mathbb{R}^d)\to H^m (\mathbb{R}^d)}\le C\varepsilon
,
\end{equation}
where the constant $C$ depends on the spatial dimension and the constants in (\ref{2}). The operator $\mathcal{K}_\varepsilon$ is determined with the help of solutions to auxiliary problems on the periodicity cell, which are introduced to define the coefficients of the homogenized operator $\hat A$. 
We note that
\begin{equation}\label{9}
\|\varepsilon^m \mathcal{K}_\varepsilon\|_{L^2 (\mathbb{R}^d)\to H^m (\mathbb{R}^d)}{\le} c, \quad 
\| \mathcal{K}_\varepsilon\|_{L^2 (\mathbb{R}^d)\to L^2 (\mathbb{R}^d)}{\le} c.
\end{equation}
Thus, the estimate   (\ref{6}) in the operator $L^2(\rd)$-norm may be obtained from (\ref{7})  by coarsening (first, by weakening the operator norm and then transferring the term $\varepsilon^m \mathcal{K}_\varepsilon$ to the remainder, due to the second inequality in  (\ref{9})).

In \cite{P20} and \cite{PMSb},  estimates of type (\ref{7}) are used in a more delicate way, which allows to improve the $L^2$-estimate (\ref{6}) with respect to the parameter $\varepsilon$ assuming that the coefficients $a_{\alpha\beta}(y)$ in (\ref{1}) are real-valued and symmetric. In this case,
the resolvent $(\hat A+1)^{-1}$ of the homogenized operator actually approximates the resolvent $(A_\varepsilon+1)^{-1}$ 
 with a remainder of order $\varepsilon^2$ in the operator $L^2(\rd)$-norm.
   Indeed, the following estimate holds
  \begin{equation}\label{10}
\|(A_\varepsilon+1)^{-1}{-}(\hat A+1)^{-1}\|_{L^2 (\mathbb{R}^d)\to L^2 (\mathbb{R}^d)}\le C\varepsilon^2,
\end{equation}
where the constant $C$ depends on the spatial dimension and the constants in (\ref{2}).
In terms of the solutions to the problems (\ref{1}) and (\ref{5}), this result provides the estimate  
\begin{equation}\label{11}
\|u^\e-u\|_{L^2(\rd)}\le  C\varepsilon^2 \|f\|_{L^2(\rd)}
\end{equation}
(with the same constant $C$ on the right-hand side as in (\ref{10})), which specifies the convergence rate in  the long-standing result  (\ref{5a}) from 70s.
 Without the symmetry condition on coefficients, it is shown \cite{P20} that
the approximation of $\varepsilon^2$ order for the  resolvent $(A_\varepsilon+1)^{-1}$ in the $L^2(\rd)$-operator norm becomes more complicated, namely,
\begin{equation}\label{12}
(A_\e+1)^{-1}=(\hat{A}+1)^{-1}+\e\K_1+O(\varepsilon^2), 
\end{equation}
where the correcting operator $\mathcal{K}_1$ is independent of $\e$ by contrast with its counterpart in (\ref{7}).

The question arises how to construct the asymptotic expansion for $(A_\e+1)^{-1}$ similar to  (\ref{12}), i.e., with remainder term of
order $\varepsilon^2$,  but in the $(L^2(\rd){\to}H^m(\rd))$-norm. 
The answer is given in Theorem \ref{th2.2}.

{\bf 1.3.} We now outline the structure of the work. An exact formulation of our 
results on improved
$H^m$-approximations is given in Section 2 (see Theorems \ref{th2.1} and \ref{th2.2}), while Sections 3--5 are of an auxiliary or preparatory nature. A direct proof of the main result, that is, the estimate (\ref{2.10}), is presented in Section 6.

In what follows we systematically use the differentiation formula for the product
\begin{equation}\label{d1}
D^\alpha(w v)=
 \sum\limits_{\gamma\le \alpha}c_{\alpha,\gamma}D^\gamma w D^{\alpha-\gamma} v=
( D^\alpha w)  v+ \sum\limits_{\gamma< \alpha}c_{\alpha,\gamma}D^\gamma w D^{\alpha-\gamma} v
\end{equation}
for  suitably differentiable functions $v$ and $w$ with some constants  $c_{\alpha,\gamma}$, where $c_{\alpha, 0}=c_{\alpha, \alpha}=1$. The sum in  (\ref{d1}) is taken over all multiindices $\gamma$ such that $\gamma\le \alpha$ or $\gamma<\alpha$. We assume that
 $\gamma\le \alpha$ if $\gamma_i\le \alpha_i$ 
 for all $1 \le i \le d$ and  $\gamma< \alpha$ if, in addition, for at least one index $i$   we have the strict inequality $\gamma_i< \alpha_i$.
 
 Throughout the paper we 
 use the following notation. Given a 1-periodic function $b(y)$, we denote by $b^\e$ or
 $(b)^\e$
  the $\e$-periodic function of the variable $x$
obtained from $ b(y)$ by substitution $y = x/\e$, i.e.,
\begin{equation}\label{d2}
b^\e(x)=b(x/\e).
\end{equation}
For example, $N_\alpha^\e{=}N_\alpha(x/\e)$, $G^\e_{\gamma\alpha}{=}G_{\gamma\alpha}(x/\e)$, 
$(D^\b G_{\gamma\alpha})^\e{=}(D^\b G_{\gamma\alpha}(y))|_{y=x/\e}$
and so on.

\section{Improved $H^m$-approximation}\label{sec:sec2}
\setcounter{theorem}{0} 
\setcounter{equation}{0}
{\bf 2.1.}  Our purpose is to obtain an $\varepsilon^2$-order approximation of the resolvent $(A_\e+1)^{-1}$ in
 the operator $(L^2(\rd){\to}H^m(\rd))$-norm. It will be the sum of the zeroth approximation and 
 several correctors, the number of which increases in comparison with 
 the  approximation from (\ref{7}). 
  By contrast to  (\ref{7}), one should take for the zeroth approximation,  instead of $(\hat{A}+1)^{-1}$, the resolvent $(\hat{A}_\e+1)^{-1}$ of the operator $\hat{A}_\e$, which is more complicated than 
  $\hat{A}$. Namely,
\beq\label{2.1}
\hat{A}_\e=(-1)^m\sum\limits_{|\alpha|=m}D^{\alpha}(
\sum\limits_{|\beta|=m}\hat{a}_{\alpha\beta}D^{\beta}
+\e \sum\limits_{|\delta|=m+1}b_{\alpha\delta}D^{\delta}
), 
\eeq
where the constant coefficients $\hat{a}_{\alpha\beta}$, which  
are the same as in (\ref{4}), and $b_{\alpha\delta}$ are determined via 
solutions to cell problems in Section 3.
Clearly,
\[
\hat{A}_\e=\hat{A}+\e B,
\]
where the differential operator $B$ is of order $2m+1$; thereby, $\hat{A}_\e$ is obtained from the commonly used homogenized operator $\hat{A}$ by a singular perturbation.

We introduce a new version of the homogenized equation with the operator (\ref{2.1})
\beq\label{2.2}
(\hat{A}_\e+1)\hat{u}^\e=f,\quad f\in \ld.
\eeq 
Since the coefficients of the equation are constant, its solution can be obtained with the help of 
the Fourier transform. 

Applying the Fourier transform to Equation
(\ref{2.2}), we get the equality 
\beq\label{2.3}
\left(1+\Lambda(\xi)+i\e\Lambda_0(\xi)\right) F[\hat{u}^\e]=F[f],\quad i=\sqrt{-1},
\eeq
where $f(x)\to F[f](\xi)$ is the Fourier transform and
$$\Lambda(\xi){=}\sum_{|\alpha|=|\beta|=m}\hat{a}_{\alpha\beta}\xi^\alpha\xi^\beta,\quad \Lambda_0(\xi){=}
\sum_{|\alpha|=m,|\delta|=m+1} b_{\alpha\delta}\xi^\alpha\xi^\delta,\quad \xi\in \rd.$$
Given any $\xi\in \rd$ and multiindex $\alpha{=}(\alpha_1,\ldots,\alpha_d)$, we have
$\xi^\alpha=\xi_1^{\alpha_1}\ldots \xi_d^{\alpha_d}$ according to our convention (see it after (\ref{1})). 
The coefficients $\hat{a}_{\alpha\beta}$ and $b_{\alpha\delta}$ are real; therefore,
 from  (\ref{2.2})  the inequality
\beq\label{2.4}
\ild \left(1+\Lambda(\xi)\right)^2 |F[\hat{u}^\e]|^2\,d\xi\le \ild|F[f]|^2\,d\xi
\eeq
follows. By the ellipticity property of $\hat{a}_{\alpha\beta}$ inherited from the 
$a_{\alpha\beta}(y)$ (see (\ref{2}))
and the  Plancherel identity,   (\ref{2.4}) 
yields the $\e$-uniform estimate
\beq\label{2.5}
\|\hat{u}^\e\|_{H^{2m}(\rd)}\le C \|f\|_{\ld},
\eeq
where the constant $C$ depends only on the spatial dimension and the constants 
 from (\ref{2}). 

{\bf 2.2.}   We are in a position to formulate the main results of the paper. We begin with approximations of the solution to
(\ref{1}) that have the structure of two-scale expansions, which is usual in homogenization.
An approximation for the solution of   (\ref{1}) is taken in the form
\beq\label{2.6}
\tilde{u}^\e(x)=u^{,\e}(x)+\e^m\,U_m^\e(x)+\e^{m+1}\,U_{m+1}^\e(x)
\eeq
with
\beq\label{2.7}\ds{
U_m^\e(x)=\sum\limits_{|\gamma|=m} N_\gamma(x/\e)D^\gamma u^{,\e}(x),}\atop\ds{
U_{m+1}^\e(x)=\sum\limits_{|\delta|=m+1}N_\delta(x/\e)D^\delta u^{,\e}(x),}
\eeq
\beq\label{2.8}
u^{,\e}(x)=\Theta^\e\, \hat{u}^\e(x),\quad \Theta^\e=S^\e S^\e ,
\eeq
where $S^\e$  is Steklov`s smoothing,
 $N_\gamma(y)$, $N_\delta(y)$,  for all multiindices $\gamma,$ $\delta$ such that $|\alpha|=m$ and $|\delta|=m+1$, and $ \hat{u}^\e(x)$ are solutions of the problems (\ref{c2}), (\ref{c15}),
  and (\ref{5}), respectively. We recall that  the Steklov smoothing operator is defined by 
\begin{equation}\label{2.9}
(S^\e\varphi)(x)=\int_Y \varphi(x-\e\omega)\,d\omega,\quad Y=[-1/2,1/2)^d,
\end{equation}
whenever
   $\varphi\in L^1_{loc}(\rd)$.

\begin{theorem}\label{th2.1}
For the difference of the solution to (\ref{1}) and the function given in 
(\ref{2.6})--(\ref{2.8}),
the 
 estimate
\begin{equation} \label{2.10}
 \|u^\e-\tilde{u}^\e\|_{H^m(\rd)}\le C\e^2 \|f\|_{L^2(\rd)}
\end{equation}
 holds, where the constant   $C$  depends only on the dimension $ d$ and the  constants 
 from (\ref{2}). 
\end{theorem}
The proof is in Section 6.

Together with $\tilde{u}^\e$, we consider the following approximation of a simpler structure:
\beq\label{2.11}
v^\e(x){=}\hat{u}^\e
{+}\e^m\,\sum\limits_{|\gamma|=m} N_\gamma(\f{x}{\e}) D^\gamma
S^\e\hat{u}^\e(x){+}\e^{m+1}
\sum\limits_{|\delta|=m+1}N_\delta(\f{x}{\e}) D^\delta S^\e S^\e\hat{u}^\e(x),
\eeq
where 
the zeroth  approximation does not involve smoothing operators, and only the first and the second correctors contain
the Steklov smoothing operator $S^\e$ defined in (\ref{2.9}) or its iteration $S^\e S^\e$ respectively.
Using only  the properties of smoothing (see Section 4), we derive  from Theorem \ref{th2.1}  the estimate
\begin{equation} \label{2.12}
 \|u^\e-{v}^\e\|_{H^m(\rd)}\le C\e^2 \|f\|_{L^2(\rd)},
\end{equation}
where the constant $C$ is of the same type as in (\ref{2.10}).

It is possible to rewrite estimates (\ref{2.10}) and (\ref{2.12}) in the operator terms, i.e., for the difference between the resolvent $(A_\e+1)^{-1}$ and its approximations. For example,
from (\ref{2.11}) and (\ref{2.12}), we deduce
\begin{equation}\label{2.13}
\|(\hat{A}_\e+1)^{-1}+\e^m\K_2(\varepsilon)+\e^{m{+}1}\K_3(\varepsilon)-(A_\e+1)^{-1}\|_{L^2 (\mathbb{R}^d)\to H^m (\mathbb{R}^d)}\le C \e^2, 
\end{equation}
where the structure of $\K_2(\varepsilon)$ and $\K_3(\varepsilon)$ is restored from  the view of the correctors in (\ref{2.11}): 
\begin{equation}\label{2.14}
\ds{
\K_2(\varepsilon)f(x)=\sum\limits_{|\gamma|=m} N_\gamma(\f{x}{\e})S^\e D^\gamma \hat{u}^{\e}(x),\, \hat{u}^{\e}(x)=(\hat{A}_\e{+}1)^{-1}f(x),
}
\atop\ds{
\K_3(\varepsilon)f(x){=}\!\sum\limits_{|\delta|=m+1}N_\delta(\f{x}{\e}) D^\delta S^\e S^\e\hat{u}^\e(x),
\, \hat{u}^{\e}(x){=}(\hat{A}_\e+1)^{-1}f(x),
} 
\end{equation}
\begin{theorem}\label{th2.2}
Let  $\hat{A}_\e$, $\K_2(\varepsilon)$, and $\K_3(\varepsilon)$ be  operators defined in (\ref{2.1}) and (\ref{2.14}). Then the sum 
 \begin{equation}\label{2.15}
 (\hat{A}_\e+1)^{-1}+\e^m\K_2(\varepsilon)+\e^{m{+}1}\K_3(\varepsilon)
 \end{equation} 
approximates the resolvent $(A_\e+1)^{-1}$ of the operator in (\ref{1}) in the energy operator norm with the 
estimate (\ref{2.13}).
\end{theorem}

{\bf 2.3.} A few words about the statement  of the above theorems and the technique we use.
It is not at once apparent that the functions given in (\ref{2.6}) and (\ref{2.11}) belong to the Sobolev space
$H^m(\rd)$. More exactly, it is not clear why the correctors in (\ref{2.6}) and (\ref{2.11}), together with their gradients of order up to $m$, belong to $\ld$ 
at first sight. But this is true thanks to the properties of the smoothing operators included inside the correctors.
Similarly, it is not quite evident that the correcting operators $\K_2(\varepsilon)$ and $\K_3(\varepsilon)$ send functions $f\in \ld$ to elements of $H^m(\rd)$. In fact, it holds that
   \begin{equation}\label{2.16}
\|\varepsilon^m \mathcal{K}_2(\e)\|_{L^2 (\mathbb{R}^d)\to H^m (\mathbb{R}^d)}{\le} C, \quad 
\|\varepsilon^{m} \mathcal{K}_3(\e)\|_{L^2 (\mathbb{R}^d)\to H^m (\mathbb{R}^d)}{\le} C
\end{equation}
with the constant $C$ depending only on the spatial dimension and the  constants in (\ref{2}). 
Being the counterparts of (\ref{9})$_1$, these  inequalities 
show that, generally,
the correctors $\e^{m}\K_2(\varepsilon)$ and $\e^{m{+}1}\K_3(\varepsilon)$  may not be transferred to the remainder term,
 for its operator norms are  of irrelevant order with respect to $\e$. Besides, by (\ref{2.16}), the orders of smallness for the correctors in (\ref{2.15}) are not alike.
 
 We pay attention to the fact  that the iterated Steklov smoothing operator is employed in (\ref{2.7}) and (\ref{2.14}); moreover, one cannot rule it out  in these formulas at all and treat only the operator $S^\e$ everywhere. 
  The necessity of the  Steklov smoothing   operator $S^\e$ and its iterations will be clarified in Section 4.

\section{Cell problems}\label{sec:sec3}
\setcounter{theorem}{0} 
\setcounter{equation}{0}

{\bf 3.1.}	On the set of smooth 1-periodic functions
$u{\in}C_\per^\infty(Y)$ with zero mean
\[
\langle u\rangle=\ilb u(y)\,dy=0,
\]  
the norm is determined by the expression
\[
\left(\ilb \sum_{|\alpha|=m}|D^{\alpha}u|^2\,dy\right)^{1/2}.\]
  The completion of this set in this norm is denoted by $\mathcal{W}$.

Since the coefficients $\{a_{\alpha\beta}(x)\}$ are periodic, then  the inequality (\ref{2}) on smooth compactly
supported functions yields a similar inequality on periodic functions:
  \beq\label{c1}
\ilb\sum\limits_{|\alpha|=|\beta|=m}a_{\alpha\beta}(y)D^{\beta}u
 D^{\alpha}u\,dy\ge  \lambda_0
\ilb \sum\limits_{|\alpha|=m}| D^{\alpha}u|^2\,dy\quad \forall u\in C_\per^\infty(Y),
\eeq
(see Lemma 3.1 in \cite{P16}).
The inequality (\ref{c1}) can be extended by closure to    $u{\in}\mathcal{W}$; thereby,
 the operator
\[
A=(-1)^m\sum_{|\alpha|=|\beta|=m}D^{\alpha}(a_{\alpha\beta}(y)D^{\beta})
\] 
acting from  $\mathcal{W}$ to its dual $ \mathcal{W}^\prime$ is coercive.

For any multiindex $\gamma$,  $|\gamma|{=}m$, we consider the problem on the cell
\beq\label{c2}
N_\gamma\in \mathcal{W}, \quad
\sum_{|\alpha|=|\beta|=m}D^{\alpha}(a_{\alpha\beta}(y)D^{\beta}N_\gamma(y))=-\sum_{|\alpha|=m}D^{\alpha}(a_{\alpha\gamma}(y)).
\eeq  
The right-hand side of (\ref{c2})  naturally defines the functional $F_\gamma$
on $\mathcal{W}$, and the equation can be written as 
\beq\label{c2a}
AN_\gamma{=}F_\gamma \quad (N_\gamma{\in} \mathcal{W}).
\eeq
Therefore, the unique solvability of (\ref{c2}) with the estimate for the solution 
\beq\label{c2b}
\|N_\gamma\|_{\mathcal{W}}\le c, \quad c=const(\lambda_0,\lambda_1),
\eeq 
is guaranteed by the following well known abstract assertion for operators acting from a real reflexive Banach space $V$ to its dual $V^\prime$.

\begin{theorem}\label{th3.1}  
Let $L:V\to V^\prime$ be a continuous linear operator satisfying the coercitivity condition:
$\langle Lv,v\rangle\ge \lambda\|v\|^2_V$
for any $v{\in}V$. 
 Then the equation $Lv{=}f\,\, (v{\in}V)$
is uniquely solvable for any $f{\in}V^\prime$,
and the solution satisfies the estimate
$
 \|v\|_V \le\lambda^{-1}\|f\|_{V^\prime}.
$
\end{theorem}

The coefficients of the homogenized operator $\hat{A}$ defined in
 (\ref{4})  are found through the solutions to the problems (\ref{c2}) as follows:
\begin{equation}\label{c3}
\hat{a}_{\alpha\beta}=
\langle a_{\alpha\beta}(\cdot)+\sum\limits_{|\gamma|=m}
a_{\alpha\gamma}(\cdot)D^{\gamma}N_\beta(\cdot)\rangle,\quad |\alpha|= m,|\beta|= m.
\end{equation}
We set 
$$
e_{\alpha\beta}
=\left\{
\begin{array}{rcl}
1,\text{ if }\alpha=\beta,\\
0,\text{ if }\alpha\neq\beta .\\
\end{array}\right.
$$
Then
\begin{equation}\label{c4}
\hat{a}_{\alpha\beta}=
\langle \sum\limits_{|\gamma|=m}
a_{\alpha\gamma}(\cdot)(e_{\gamma\beta}+D^{\gamma}N_\beta(\cdot))\rangle, 
\end{equation}
or 
\begin{equation}\label{c5}
\hat{a}_{\alpha\beta}=\langle \tilde{a}_{\alpha\beta}\rangle,\quad \tilde{a}_{\alpha\beta}(y)=
 \sum\limits_{|\gamma|=m}
a_{\alpha\gamma}(y)(e_{\gamma\beta}+D^{\gamma}N_\beta(y)).
\end{equation}

It is known (see, e.g., [3, Lemma 3.2]) that the matrix  of the homogenized coefficients (\ref{c3}) belongs to the class (\ref{2}). This fact is essential in the proof of the elliptic estimate (\ref{2.5}).

Setting 
 \beq\label{c7}
g_{\alpha\beta}(y)=\tilde{a}_{\alpha\beta}(y)-\hat{a}_{\alpha\beta} \quad \forall \a,\b,
\eeq
we obtain the relations
\beq\label{c8}
 \langle g_{\alpha\beta}\rangle=0,\quad 
\sum\limits_{|\a|=m}D^\a g_{\alpha\beta}=0 \quad \forall \b, 
\eeq
which follow from the statement of the problem (\ref{c2}) and the definitions (\ref{c7}) and (\ref{c5}).
  The equality  (\ref{c8})$_2$ can be understood in the sense of the integral identity on periodic functions
\beq\label{c9}
\sum\limits_{|\a|=m}\ilb  g_{\alpha\beta} D^\a\varphi\,dy=0 \quad\forall\varphi\in C^\infty_\per(Y),
\eeq
or in the sense of distributions in $\rd$. The integral identity (\ref{c9}) is extended by closure to functions $\varphi\in \mathcal{W}$.

{\bf 3.2.} We recall that periodic solenoidal vectors with zero mean admit representations in terms of the divergence of a skew-symmetric matrix. More exactly, if
$g{\in} L^2_\per(Y)^d$,  $\div\, g{=}0$, and $\langle g\rangle{=}0$,  then there exists a skew-symmetric matrix
 $G\in H^1_\per(Y)^{d\times d}$ such that 
 $$\div \, G{=}g, \quad \|G\|_
{H^1_\per(Y)^{d\times d}}\le c\|g\|_{H^1_\per(Y)^d} , \,c{=}const(d)$$ (see, e.g., \cite{JKO}, Chapter 1). The property $\div\, g{=}0$ is understood in the sense of the integral identity
\beq\label{c10}
\langle g\cdot \nab \varphi\rangle=0 \quad \forall \varphi\in C^\infty_\per(Y),
\eeq
or in the sense of distributions in $\rd$.

The following lemma extends 
the above assertion on solenoidal vectors (in the classical sense) to the case where the orthogonality property of type (\ref{c10}) holds with the gradient $\nab^m$ of order $m\ge 2$; here, $\nab^m \varphi=\{D^\a\varphi\}_{\a},$ $ |\a|{=}m$.
Let $p$ be the number of multiindices of length $m$, then $\nab^m\varphi$ is regarded as a vector with $p$ components.  

\begin{lemma}\label{lem2.1} Assume that $\{g_\alpha\}_{|\alpha|=m}
\in L_\per^2(Y)^p$ 
and
\beq\label{c11}
\langle g_\alpha\rangle=0,\quad
\sum_{|\a|=m}D^\a g_\a=0.
\eeq
Then there exists a matrix $\{G_{\gamma\alpha}\}_{|\alpha|=|\gamma|=m}$ from $H^m_\per(Y)^{p\times p}$ such that
for all miltiindices $\alpha$, $\gamma$
\beq\label{c12}
G_{\gamma\alpha}=-G_{\gamma\alpha}, \quad \sum_{|\gamma|=m}D^\gamma G_{\gamma\alpha}=g_\alpha,
\eeq
\beq\label{c13}
\|G_{\gamma\alpha}\|_{H^m(Y)}\le c \sum_{|\alpha|=m}\|g_\alpha\|_{L^2(Y)},\quad c=const( {d,m}).
\eeq
\end{lemma}
 The proof is given in \cite{P16} (see also   \cite{AA16} and \cite{UMN}). 
 
 For every fixed $\beta$ the vector
 $\{g_{\alpha\beta}\}_{\a}$, $|\a|=m$, in (\ref{c7}) satisfies the assumptions of Lemma \ref{lem2.1}. Consequently, there is a matrix
 $\{G_{\gamma\alpha\beta}\}_{\gamma,\alpha}\in H^m$, ${|\a|{=}|\gamma|{=}m}$,  such that identities of the form (\ref{c12}) hold 
 componentwise, i.e.,
 \beq\label{c14}
G_{\alpha\gamma\beta}=-G_{\gamma\alpha\beta}, \quad
g_{\alpha\beta}=\sum_{|\gamma|=m}D^\gamma G_{\gamma\alpha\beta},
\eeq
and $ G_{\gamma\alpha\beta}$  satisfies $H^m$-estimate of kind (\ref{c13}). 

We are ready now to introduce, for any multiindex $\delta$,  $|\delta|{=}m{+}1$,  the problem on the cell, which depends on the solutions $N_\gamma $ to (\ref{c2}) and the functions $ G_{\gamma\alpha\beta}$ from the
 representation (\ref{c14}). This  problem is a kind of (\ref{c2a}), namely,
\beq\label{c15}
N_\delta\in \mathcal{W}, \quad
\sum_{|\alpha|=|\beta|=m}D^{\alpha}(a_{\alpha\beta}D^{\beta}N_\delta(y))=-\sum_{|\alpha|=m}D^{\alpha}
F_{\alpha,\delta}.
\eeq  
The right-hand side in (\ref{c15})  
defines the functional
on $\mathcal{W}$ if we set 
\beq\label{c16}
F_{\alpha,\delta}=\sum_{|\gamma|=|\beta|=m,\,\b<\delta}c_{\gamma,\beta+\gamma-\delta}(a_{\alpha\gamma}(y)D^{\beta+\gamma-\delta}N_\b+D^{\beta+\gamma-\delta}G_{\a\gamma\b})
\eeq
with the constants $c_{\beta,\beta+\gamma-\delta}$ from the product rule (\ref{d1}).

Related to (\ref{c15}) and (\ref{c16}), the matrix
 \beq\label{c17}
\tilde{g}_{\alpha\delta}(y)=
\sum_{|\beta|=m}a_{\alpha\beta}(y)D^{\beta}N_\delta(y)+F_{\alpha,\delta}-b_{\alpha\delta},
\eeq
where
 \beq\label{c18}
b_{\alpha\delta}=\langle 
\sum_{|\beta|=m}a_{\alpha\beta}(y)D^{\beta}N_\delta(y)+F_{\alpha,\delta}\rangle,
\eeq
turns to be the counterpart of $g_{\alpha\beta}(y)$. In fact, 
for every fixed $\delta$, $|\delta|{=}m{+}1$, the vector
 $\{\tilde{g}_{\alpha\delta}\}_{\a}$, $|\a|{=}m$, in (\ref{c17}) satisfies the assumptions of Lemma \ref{lem2.1},
 namely,
 \beq\label{c19a}
 \langle \tilde{g}_{\alpha\delta}\rangle=0,\quad\sum_{|\a|=m}D^\a  \tilde{g}_{\alpha\delta}=0.
 \eeq
  Consequently, there is a matrix
 $\{\tilde{G}_{\gamma\alpha\delta}\}_{\gamma,\alpha}$, ${|\a|{=}|\gamma|{=}m}$,  such that 
 \beq\label{c19}
\tilde{G}_{\alpha\gamma\delta}=-\tilde{G}_{\gamma\alpha\delta}, \quad
\tilde{g}_{\alpha\delta}=\sum_{|\gamma|=m}D^\gamma \tilde{G}_{\gamma\alpha\delta},
\eeq
and $ \tilde{G}_{\gamma\alpha\delta}$  satisfies $H^m$-estimate of kind (\ref{c13}). 

\begin{rem}\label{remark4.1} \rm
Taking into account the structure of the functions $F_{\alpha,\delta}$ defined in (\ref{c16}), we can simplify the formula for the coefficients $b_{\alpha\delta}$ in such a way that  $b_{\alpha\delta}$ will be expressed only in terms of the solutions to  
(\ref{c2}) and (\ref{c15}).
\end{rem}

\begin{rem}\label{remark4.2} \rm 
   We can rewrite equation (\ref{c15}) for $N_\delta$, $|\delta|=m+1$, using firstly (\ref{c14})$_2$, then (\ref{c7}) and (\ref{c5}), thus, excluding in the right-hand side  firstly the functions  $G_{\a\gamma\b}$ and then $g_{\a\b}$.
As a result,  the right-hand side  in  (\ref{c15}) will be determined only by the functions $N_\gamma$, 
$|\gamma|=m$, 
 its derivatives,  and the coefficients $a_{\a\b}$. For further calculations,
 we prefer the equation for $N_\delta$, $|\delta|=m+1$, in display (\ref{c15}).
\end{rem}

{\bf 3.3.} 
 As a corollary of  Lemma \ref{lem2.1}, we have the next 
\begin{lemma}\label{lem2.2}
Let a 1-periodic vector $\{g_\alpha\}_{|\alpha|=m}
\in L_\per^2(Y)^p$ 
satisfy    the assumptions of Lemma \ref{lem2.1}, and let a 1-periodic matrix 
$ \{G_{\gamma\alpha}(y)\}_{|\alpha|=|\gamma|=m}$  be the associated matrix potential such that  (\ref{c12})  and (\ref{c13}) hold.
Then 
\begin{equation}\label{c20}
\ds{g_\alpha(x/\e)\Phi(x)}
\atop\ds{
=\sum_{|\gamma|=m} D^\gamma (\e^m G^\e_{\gamma\alpha} \Phi)
-\sum_{|\gamma|=m}\sum_{\mu<\gamma}\e^{m-|\mu|}c_{\gamma,\mu} (D^\mu  G_{\gamma\alpha})^\e D^{\gamma-\mu}\Phi}
\end{equation}
for any
$\Phi \in \C0$  and all indices $\alpha$, $|\alpha|=m$, where the constants $c_{\gamma,\mu}$
are from (\ref{d1}). Furthermore, the vector
 \begin{equation} \label{c21}
\{M_\alpha(x)\}_{|\alpha|=m},\quad
M_\alpha(x)=\sum_{|\gamma|=m} D^\gamma (G^\e_{\gamma\alpha} \Phi),
\end{equation}
satisfies the relation of kind (\ref{c11})$_2$, that is, 
 \begin{equation} \label{c22}
\sum_{|\a|=m} D^\a M_\a=0 \quad (\mbox{in the sense of distributions on}\quad\rd).
\end{equation}
\end{lemma}
\textbf{Proof.}
By assumption,
\[
g_\alpha(y)
=\sum_{|\gamma|=m} D^\gamma G_{\gamma\alpha}(y),\quad
g_\alpha(x/\e)\Phi(x)=\sum_{|\gamma|=m} D^\gamma (\e^m G_{\gamma\alpha}(x/\e)) \Phi(x).
\]
Hence, to obtain (\ref{c20}), it suffices to recall the product rule (\ref{d1}) taking into account the two-scale nature of the expression $g_\alpha(x/\e)\Phi(x)$ and using the notation  (\ref{d2}).

The relation (\ref{c22}) holds for the vector (\ref{c21})
because of the integral identity
\[
\ild \varphi \sum_{|\a|=m} D^\a M_\a\,dx=0\quad\forall \varphi\in \C0
\]
obtained from the chain of integral identities
\[
\ild \varphi \sum_{|\a|=m} D^\a M_\a\,dx=
\ild \varphi \sum_{|\gamma|=|\a|=m}D^\gamma D^\a (G^\e_{\gamma\a}\Phi)\,dx
\]
\[=
\ild \sum_{|\gamma|=|\a|=m}(D^\gamma D^\a  \varphi )G^\e_{\gamma\a}\Phi\,dx
=\ild \Phi\sum_{|\gamma|=|\a|=m}(D^\gamma D^\a  \varphi )G^\e_{\gamma\a}\,dx,
\]
and the pointwise identity 
\[
\sum_{|\gamma|=|\a|=m}(D^\gamma D^\a  \varphi )G^\e_{\gamma\a}=0
\]
since  the matrix $G_{\gamma\a}$ is skew-symmetric (see (\ref{c12})$_1$).

\section{Smoothing operators}\label{sec:sec4}
\setcounter{theorem}{0} 
\setcounter{equation}{0}

 In this paper, we prove homogenization estimates by method, coming from \cite{Zh1} and \cite{ZhP05} (see also the 
overview \cite{UMN}), where  the difficulties connected with the lack of regularity in data are overcome by using smoothing operators, among which we  highlight out Steklov`s smoothing operator $S^\e$ defined in (\ref{2.9}).
We start with  
the simplest 
properties of the operator $S^\e$:
  \begin{equation}\label{m.2}
\|S^\e\varphi\|\le\|\varphi\|,
\end{equation}
 \begin{equation}\label{m.3}
\|S^\e\varphi-\varphi\|\le (\sqrt{d}/2)\e\|\nab\varphi\|\quad \forall\varphi\in H^1(\R^d).
\end{equation}
Here  and throughout this Section, $\|\cdot\|$ denotes, for brevity,  the  norm   in $\ld$.
We mention also the evident property $S^\e(D^\a \varphi)=D^\a 
(S^\e\varphi)$ for any derivative $D^\a$, which is exploited systematically in the sequel.

The following property of the Steklov smoothing operator $S^\e$ to interrelate with $\e$-periodic factors  is pivotal for our method (see its simple proof, e.g., in \cite{ZhP05} or
\cite{UMN}). 
  \begin{lemma}\label{LemM1} If $\varphi{\in}L^2(\R^d)$, $b{\in}L^2_\per(Y)$,
 and $b^\e(x){=}b(\e^{-1}x)$,
 then 
 \begin{equation}\label{m.4}
\|b^\e S^\e\varphi\|\le\langle b^2\rangle^{1/2}\|\varphi\|.
\end{equation}
\end{lemma} 

The above estimate (\ref{m.3}) 
can be improved under higher regularity conditions. For example,  
 \begin{equation}\label{m.6}
\|S^\e\varphi-\varphi\|\le C\e^2\|\nab^2\varphi\|\quad \forall\varphi\in H^2(\R^d),\quad C=const(d).
\end{equation}
By duality, from (\ref{m.6}) we obtain 
\begin{equation}\label{m.7}
\|S^\e\varphi-\varphi\|_{H^{-2}(\rd)}\le C\e^2\|\varphi\|_{L^2(\rd)}\quad\forall\varphi\in L^2(\R^d),\quad C=const(d).
\end{equation}

 Instead of Steklov`s smoothing operator, we can often consider the general  smoothing operator with  an arbitrary kernel. This is
 \begin{equation}\label{m.11}
\Theta^\e\varphi(x)=\ild\varphi(x-\e\omega)\theta(\omega)\,d\omega, 
\end{equation}
where $\theta\in L^\infty(\rd)$ 
is compactly supported, $\theta\ge 0$, and $\int_{\rd} \theta(x)dx{=}1$. 

The estimates  (\ref{m.2})--(\ref{m.4}) formulated for the Steklov   smoothing operator remain valid for the general smoothig operator
 (\ref{m.11}) with only one note that right-hand side contains constants depending not only on the dimension $d$, but also on the kernel $\theta$. If, in addition,  the kernel $\theta$ is even, 
then $\Theta^\e$ possesses the properties of  type (\ref{m.6}) and (\ref{m.7}).

The following useful properties of the operator (\ref{m.11}) are highlighted in  \cite{NY}  and  \cite{P20s}. 
\begin{lemma}\label{LemM2} Assume that $\theta$ is piecewise $C^k$-smooth, $k$ is a natural number, 
$b\in L^2_\per(Y)$, $b^\e(x)=b(x/\e)$, and  
  $\varphi\in L^2(\rd)$. Then
   \begin{equation}\label{m.12}
\|\Theta^\e\nab^k\varphi\|\le C\e^{-k}
\|\varphi\|,\quad C=const(\theta),
\end{equation}
   \begin{equation}\label{m.13}
\|b^\e\Theta^\e\nab^k\varphi\|\le C\e^{-k}
\langle b^2\rangle^{1/2}\|\varphi\|,\quad C=const(\theta).
\end{equation}
\end{lemma}

It is obvious that the Steklov smoothing operator $S^\e$ can be defined by (\ref{m.11}) if the smoothing kernel  $\theta_1(x)$ is the characteristic function of the cube   $Y{=}[-{1}/{2},{1}/{2})^d$.
The double Steklov smoothing operator $(S^\e)^2=S^\e S^\e$ has the form (\ref{m.11}) with the smoothing kernel equal to the convolution  $\theta_2=\theta_1*\theta_1$.
The kernels $\theta_2$  (see computations for it in \cite{P20s}) turns to be piecewise $C^1$-smooth. As a consequence, the properties  (\ref{m.12}) and (\ref{m.13})  with $k{=}1$ hold 
for the smoothing operator $\Theta^\e=(S^\e)^2$. In addition, since $\theta_2$  is even,  the operator   $(S^\e)^2$  possesses the properties of type (\ref{m.6}) and (\ref{m.7}).

\section{Preliminary calculations}
\setcounter{theorem}{0} 
\setcounter{equation}{0}

Before we proceed to the direct proof of the estimate (\ref{2.10}), we make some 
formal calculations under strong regularity assumptions which will be later taken off. Let
\begin{equation}\label{5.1}
v^\e(x)=u(x)+\e^m
 \sum\limits_{|\gamma|= m}N^\e_\gamma(x)D^\gamma u(x)+ 
 \e^{m+1}
 \sum\limits_{|\delta|= m+1}N_\delta^\e(x)D^\delta u(x),
\end{equation}
where  $N_\gamma$, $|\gamma|= m$, and $N_\delta$, $|\delta|= m+1$, are the solutions to the cell problems (\ref{c2})
and (\ref{c15}) respectively.

 We first suppose that the function $u(x)$ in (\ref{5.1}) is infinitely differentiable and, together with its derivatives, is decreasing at infinity sufficiently rapidly, so that 
$v^\e\in H^m(\rd)$. Therefore,   the discrepancy of the function
$v^\e$ in (\ref{1}), that is $(A_\e+1)v^\e-f$, can be calculated. But let 
the right-hand side function $f$  
be of the form $f=(\hat{A}_\e+1)u$, where $\hat{A}_\e$ is the homogenized operator defined in (\ref{2.1}) and the function $u(x)$ is the same as in (\ref{5.1}). Thus,
\begin{equation}\label{5.2}\ds{
(A_\e+1)v^\e-f=(A_\e+1)v^\e-(\hat{A}_\e+1)u=A_\e v^\e-\hat{A}_\e u +(v^\e-u)}
\atop\ds{\stackrel{(\ref{1}),(\ref{2.1})}=
\sum_{|\a|= m}(-1)^m D^\a(\Gamma_\a(v^\e,A_\e)-\Gamma_\a(u,\hat{A}_\e))+(v^\e-u),
}
\end{equation}
where we introduce 
 the generalized gradients
 \beq\label{5.3}\ds{
\Gamma_\a(v^\e,A_\e)=
 \sum\limits_{|\beta|=  m}a^\e_{\alpha\beta}D^{\beta}v^\e,}
 \atop\ds{
 \Gamma_\a(u,\hat{A}_\e){=} \Gamma_\a(u,\hat{A}){+} \Gamma_\a(u,\e B){=}
\! \sum\limits_{|\beta|= m}\hat{a}_{\alpha\beta} D^{\beta}u{+}
\! \sum\limits_{|\delta|= m+1}\e b_{\alpha\delta} D^{\delta}u}
 \eeq 
 for any index $\a$, $|\a|= m$.

 In view of (\ref{5.2}), we need to compare the generalized gradients.
 A simple calculation shows that, for any $\b$, $|\b|=m,$
\[
D^{\beta}(\e^m N^\e_\gamma D^{\gamma}u)=(D^{\beta}N_\gamma)^\e D^{\gamma}u+
\sum_{\mu<\b}\e^{m-|\mu|}c_{\b,\mu} (D^\mu  N_\gamma)^\e D^{\b+\gamma-\mu}u,
\]
by the product rule (\ref{d1}). Thus, due to (\ref{5.3})$_1$ and 
(\ref{5.1}), we have
\[
\Gamma_\a(v^\e,A_\e)
=
 \sum_{|\beta|=  m}a^\e_{\alpha\beta}D^{\beta}
( u+
\e^m \!\!\sum_{|\gamma|= m}N_\gamma^\e D^\gamma u+\e^{m+1}
 \!\!\sum\limits_{|\delta|{=} m{+}1}N_\delta^\e D^\delta u
)\]
\[=
\!\sum_{|\beta|= m}
( a^\e_{\alpha\beta}D^{\beta}u{+}
 \sum_{|\gamma|=m}a^\e_{\alpha\gamma}(D^{\gamma}N_\beta)^\e D^{\beta}u)+
 \e\sum_{|\gamma|=m}\sum\limits_{|\delta|{=} m{+}1}\!a^\e_{\alpha\gamma}(D^{\gamma}N_\delta)^\e D^{\delta}u
  \]
\[+\e\sum_{|\gamma|=|\beta|= m}\, \sum_{\mu<\gamma,|\mu|=m-1}a^\e_{\alpha\gamma}c_{\gamma,\mu}(D^{\mu}N_\beta)^\e D^{\beta+\gamma-\mu}u+w^\e_\a,
\]
where the sum
\begin{equation}\label{5.4}
\ds{
w^\e_\a=
\sum_{|\gamma|=|\beta|= m}\, \sum_{\mu<\gamma,|\mu|<m-1}\e^{m-|\mu|}a^\e_{\alpha\gamma}c_{\gamma,\mu}(D^{\mu}N_\beta)^\e D^{\beta+\gamma-\mu}u}
\atop\ds{
+\sum_{|\gamma|=m}\sum\limits_{|\delta|{=} m{+}1}\sum_{\mu<\gamma,|\mu|<m}\e^{m+1-|\mu|}
a^\e_{\alpha\gamma}c_{\gamma,\mu}(D^{\mu}N_\delta)^\e D^{\delta+\gamma-\mu}u
}
\end{equation}
collects all the terms containing  powers  $\e^n$, $n\ge 2$, as factors.

Using the notation from (\ref{c4})--(\ref{c7}), we  transform the first sum in the above representation of $\Gamma_\a(v^\e,A_\e)$ as follows:
\[
\sum_{|\beta|= m}
( a^\e_{\alpha\beta}D^{\beta}u{+}
 \sum_{|\gamma|=m}a^\e_{\alpha\gamma}(D^{\gamma}N_\beta)^\e D^{\beta}u)=
 \! \sum_{|\beta|= |\gamma|=m}a^\e_{\alpha\gamma}(e_{\gamma\beta}+(D^{\gamma}N_\beta)^\e) D^{\beta}u
\]
\[
\stackrel{ (\ref{c4}),(\ref{5.3})_2}=\sum_{|\beta|=  m}
(\tilde{a}^\e_{\alpha\beta}-\hat{a}_{\alpha\beta})D^{\beta}u+\Gamma_\a(u,\hat{A})=
\Gamma_\a(u,\hat{A})+
 \sum_{|\beta|=  m}
g^\e_{\alpha\beta}
D^\beta u.
\]
Consequently, we rewrite the representation 
\[
\Gamma_\a(v^\e,A_\e)=\Gamma_\a(u,\hat{A})+
 \sum_{|\beta|=  m}g^\e_{\alpha\beta}
D^\beta u
+ \e\sum_{|\gamma|=m}\sum\limits_{|\delta|{=} m{+}1}a^\e_{\alpha\gamma}(D^{\gamma}N_\delta)^\e D^{\delta}u\]
\begin{equation}\label{5.5}
+\e\sum_{|\gamma|=|\beta|= m}\, \sum_{\mu<\gamma,|\mu|=m-1}a^\e_{\alpha\gamma}c_{\gamma,\mu}(D^{\mu}N_\beta)^\e D^{\beta+\gamma-\mu}u+w^\e_\a.
\end{equation}

Applying Lemma \ref{lem2.2} to the term $g^\e_{\alpha\beta}
D^\beta u$ (note that the vector $\{g_{\alpha\beta}\}_\a$, with $\b$ fixed, satisfies the assumptions of Lemma \ref{lem2.2}),
we  obtain 
\begin{equation}\label{5.6}
\ds{
g^\e_{\alpha\beta}
D^\beta u=\sum_{|\gamma|=  m}(D^\gamma G_{\gamma\a\b})^\e D^\beta u}\atop\ds{
=\sum_{|\gamma|=m} D^\gamma (\e^m G^\e_{\gamma\alpha\b} D^\beta u)
-\sum_{|\gamma|=m}\sum_{\mu<\gamma}\e^{m-|\mu|}c_{\gamma,\mu} (D^\mu  G_{\gamma\alpha\b})^\e D^{\gamma-\mu}D^\beta u.
}
\end{equation}
Setting
\[
M_\a:=\sum_{|\gamma|=m} D^\gamma ( G^\e_{\gamma\alpha\b} D^\beta u),
\]
we get the vector $\{M_\alpha(x)\}_{|\alpha|=m}$ with the property (\ref{c22}).
Therefore, 
\begin{equation}\label{5.7}
\ds{
\sum_{|\a|= m} D^\a(g^\e_{\alpha\beta}
D^\beta u)\stackrel{(\ref{5.6})}=
{-}\!\sum_{|\a|=|\gamma|= m}D^\a
\sum_{\mu<\gamma}\e^{m-|\mu|}c_{\gamma,\mu} (D^\mu  G_{\gamma\alpha\b})^\e D^{\beta+\gamma-\mu}u
}\atop\ds{
=\!\sum_{|\a|=|\gamma|= m}D^\a
\,
\sum_{\mu<\gamma, |\mu|=m-1}\e c_{\gamma,\mu} (D^\mu  G_{\alpha\gamma\b})^\e D^{\beta+\gamma-\mu} u
+\sum_{|\a|= m}D^\a
{w}^\e_{\a\b},}
\end{equation}
where the sum
\begin{equation}\label{5.8}
{w}^\e_{\a\b}=\sum_{|\gamma|=m}\sum_{\mu<\gamma, |\mu|<m-1}\e^{m-|\mu|} c_{\gamma,\mu} (D^\mu  G_{\alpha\gamma\b})^\e D^{\beta+\gamma-\mu}u
\end{equation}
collects the terms with powers $\e^n$, $n\ge 2$, as a factor.
To change the sign in (\ref{5.7}), the skew-symmetry property of $ G_{\alpha\gamma\b}$ is used, i.e., $ G_{\alpha\gamma\b}=-G_{\gamma\alpha\b}$.

 From (\ref{5.3}),  
 (\ref{5.5}), and (\ref{5.7}) it follows that
\begin{equation}\label{5.9}
R_\e:=\sum_{|\a|= m} D^\a(\Gamma_\a(v^\e,A_\e)-\Gamma_\a(u,\hat{A}_\e))
\end{equation}
\[
=\sum_{|\a|= m} D^\a(\Gamma_\a(v^\e,A_\e)-\Gamma_\a(u,\hat{A})-\Gamma_\a(u,\e B))
\]
\[
=\e\sum_{|\a|=|\gamma|= m}D^\a \sum\limits_{|\delta|{=} m{+}1}a^\e_{\alpha\gamma}(D^{\gamma}N_\delta)^\e D^{\delta}u\]
\[
+\e\sum_{|\a|=|\gamma|=|\beta|= m}D^\a( \sum_{\mu<\gamma,|\mu|=m-1}a^\e_{\alpha\gamma}c_{\gamma,\mu}(D^{\mu}N_\beta)^\e D^{\beta+\gamma-\mu}u
\]
\[
+\sum_{\mu<\gamma, |\mu|=m-1} c_{\gamma,\mu} (D^\mu  G_{\alpha\gamma\b})^\e D^{\beta+\gamma-\mu} u)
-\e\sum\limits_{|\a|= m}D^\a
\sum\limits_{|\delta|= m+1} b_{\alpha\delta} D^{\delta}u
\]
\[
+\sum_{|\a|= m}D^\a w^\e_\a+\sum_{|\a|=|\b|= m}D^\a{w}^\e_{\a\b},
\]
where $w^\e_\a$ and ${w}^\e_{\a\b}$ defined in (\ref{5.4}) and (\ref{5.8}) are of order $O(\e^2)$ for small $\e$.
We can  shortly write
\[
R_\e=\e R_{\e,1}+\sum_{|\a|= m}D^\a w^\e_\a+\sum_{|\a|=|\b|= m}D^\a{w}^\e_{\a\b},
\]
where the summand $\e R_{\e,1}$ absorbs all the terms containing the factor $\e$ in the above representation.
It holds 
$D^\delta{=}D^{\beta+\gamma-\mu}$ if the exponents coincide, i.e., $\delta{=}\beta+\gamma-\mu$ or
$\mu{=}\beta+\gamma-\delta$. In this case, $\mu{<}\gamma\Leftrightarrow \b{<}\delta$. Rearranging the summation in the above representation, we get
\[
\ds{
\e R_{\e,1}=\e\sum_{|\a|=m}D^\a \sum_{\delta=m+1}D^{\delta}u
[\sum_{|\gamma|=m}(a_{\alpha\gamma} D^\gamma N_\delta
}\atop\ds{
+\sum_{|\b|=m,\b<\delta}c_{\gamma,\beta+\gamma-\delta}(a_{\alpha\gamma} D^{\beta+\gamma-\delta} N_\b+
D^{\beta+\gamma-\delta}G_{\a\gamma\b}))^\e-b_{\a\delta}
].
}
\]
The oscillating $\e$-periodic multiplier standing above in the square brackets  can be written as 
$\tilde{g}^\e_{\a \delta}$, and the corresponding 1-periodic function $\tilde{g}_{\a \delta}$ is just the same as in (\ref{c17})
according to 
(\ref{c16})--(\ref{c18}).

We summarize  in a short way what is obtained by this moment:
\begin{equation}\label{5.10}
R_\e=
\e\sum_{|\a|=m}D^\a \sum_{\delta=m+1}\tilde{g}^\e_{\a \delta}D^{\delta}u
+\sum_{|\a|= m}D^\a w^\e_\a+\sum_{|\a|=|\b|= m}D^\a{w}^\e_{\a\b}.
\end{equation}
For each fixed multiindex $\delta$, $|\delta|{=}m+1$, by properties of $\tilde{g}_{\a \delta}$
(see (\ref{c19a}) and (\ref{c19}))  
and Lemma \ref{lem2.2}, it holds
\[
\ds{\tilde{g}_{\a \delta}(x/\e)\Phi(x)}
\atop\ds{
=\sum_{|\gamma|=m} D^\gamma (\e^m\tilde{ G}^\e_{\gamma\alpha\delta} \Phi)
-\sum_{|\gamma|=m}\sum_{\mu<\gamma}\e^{m-|\mu|}c_{\gamma,\mu} (D^\mu  \tilde{ G}_{\gamma\alpha\delta})^\e D^{\gamma-\mu}\Phi}
\]
with $\Phi=D^{\delta}u$ and the matrix $\{\tilde{G}_{\gamma\alpha\delta}\}_{\gamma,\alpha}$, ${|\a|{=}|\gamma|{=}m}$, satisfying (\ref{c19}), where
\[
\sum_{|\a|=m}D^\a \sum_{|\gamma|=m} D^\gamma (\e^m\tilde{ G}^\e_{\gamma\alpha\delta} \Phi)=0.
\]
Finally, we transform 
(\ref{5.10}) as follows:
\begin{equation}\label{5.11}\ds{
R_\e=
-\e\sum_{|\a|=|\gamma|=m}D^\a \sum_{\delta=m+1}\,
\sum_{\mu<\gamma}\e^{m-|\mu|}c_{\gamma,\mu} (D^\mu  \tilde{ G}_{\gamma\alpha\delta})^\e 
D^{\gamma-\mu+\delta}u}
\atop\ds{
+\sum_{|\a|= m}D^\a w^\e_\a+\sum_{|\a|=|\b|= m}D^\a{w}^\e_{\a\b}.}
\end{equation}
Thus, each term in the right-hand side of (\ref{5.11}) includes some power $\e^n$, $n\ge 2$, as a factor
(see 
the definitions (\ref{5.4}) and (\ref{5.8}) for $w^\e_\a$ and ${w}^\e_{\a\b}$). The same is true for the the
entire right-hand side of (\ref{5.2}) if we recall the notation (\ref{5.9}) and  the equality
\[
v^\e-u\stackrel{(\ref{5.1})}=\e^m
 \sum\limits_{|\gamma|= m}N^\e_\gamma D^\gamma u+ 
 \e^{m+1}
 \sum\limits_{|\delta|= m+1}N_\delta^\e D^\delta u,
\]
where $m\ge 2$. This yields the final representation for the discrepancy of the function (\ref{5.1}) with  equation (\ref{1}) as the sum of several terms of the same type 
\begin{equation}\label{5.12}
(A_\e+1)v^\e{-}f{=}\sum_j \e^{n_j}b_j(x/\e)\Phi_j(x){+}
\!\sum_{|\a|=m}D^\a \sum_j \e^{n_j}\tilde{b}_j(x/\e)\tilde{\Phi}_j(x). 
\end{equation}
Here, any exponent $n_j\ge 2$  and all $\e$-periodic functions $b_j(x/\e)$ and $\tilde{b}_j(x/\e)$ are formed of their 1-periodic counterparts from the list
\begin{equation}\label{5.13}
a_{\a\b},\,N_\gamma,\, D^\mu N_\gamma,\, G_{\gamma\a\b},\, D^\mu G_{\gamma\a\b},\,
 \tilde{G}_{\gamma\a\b},\, D^\mu \tilde{G}_{\gamma\a\b}
\end{equation}
consisting of the coefficients of (\ref{1}), the solutions to the cell problems (\ref{c2}) and (\ref{c15}) together with the derivatives up to the order $m$, the components of the matrix potentials from (\ref{c14}) and (\ref{c19})
 together with the derivatives up to the order $m$. 
 The functions  $\Phi_j$  and $\tilde{\Phi}_j$ in (\ref{5.12}) coincide  with 
 the function $u$ or its derivatives $D^\nu u$ up to the order $2m+1$, where $(\hat{A}_\e+1)u=f$.

Since $({A}_\e{+}1)u^\e{=}f{=}(\hat{A}_\e{+}1)u$ for the solution $u^\e$  to (\ref{1}), from (\ref{5.12}) 
we  get 
\[
(A_\e+1)(v^\e-u^\e){=}O(\e^2),
\]
whereof, by the energy estimate, we  derive the estimate 
$$
\|v^\e-u^\e\|_{H^m(\rd)}=O(\e^2), 
$$
which is not the same as the desired estimate (\ref{2.10}). We remind that our preceding  calculations 
and the final estimate make sense only under strong regularity assumptions. Furthermore, 
the majorant in the extracted 
estimate cannot be guaranteed to be of the type given in (\ref{2.10}). We will show how to cope with these difficulties in the next Section.

\section{$H^m$-estimate of order $\e^2$}\label{sec:sec6}
\setcounter{theorem}{0} 
\setcounter{equation}{0}

In this Section, we will prove the  $H^m$-estimate (\ref{2.10}) for the difference between the solution $u^\e$ to  (\ref{1}) and its approximation defined in (\ref{2.6}).

{\bf 6.1.} Firstly, let us verify that the function (\ref{2.6})
belongs to the space $H^m(\rd)$. It suffices to show that every term in the both correctors (\ref{2.7}), together with
its derivatives up to order $m$, belongs to the space $\ld$. For this we apply lemmas \ref{LemM1} and \ref{LemM2}.
The latter one is employed to handle the products 
\begin{equation}\label{6.1}
\e^{m+1}N_\delta(x/\e)D^{\alpha+\delta}u^{,\e}(x),\quad |\a|=m,\, |\delta|=m+1,
\end{equation}
 that emerge from differentiation of order $m$ of the second corrector in (\ref{2.7}). In this case,
$$
N_\delta(x/\e)D^{\alpha+\delta}u^{,\e}(x)\stackrel{(\ref{2.8})}=N_\delta(x/\e)D^{\alpha+\delta}\Theta^\e\, \hat{u}^\e(x)=N_\delta(x/\e)\Theta^\e D_j\varphi(x),
$$
where $\varphi=D^\nu \hat{u}^\e$ for some multiindex $\nu$, $|\nu|=2m$, and thereby, $\varphi\in \ld$ with the estimate 
\begin{equation}\label{6.2}
\|\varphi\|_{\ld}=\|D^\nu \hat{u}^\e\|_{\ld}\stackrel{(\ref{2.5})}\le C\|f\|_{\ld}.
\end{equation}
It remains to note that the smoothing operator $\Theta^\e=S^\e S^\e$ has the smoothing kernel which is
 piecewise $C^1$-smooth (see discussion in the end of Section 4). Thus, applying (\ref{m.13}) yields
 \[
 \|N_\delta^\e D^{\alpha+\delta}u^{,\e}\|_{\ld}=\|N_\delta^\e\Theta^\e D_j\varphi\|_{\ld}\le C \e^{-1}\langle N_\delta^2\rangle^{1/2}\|\varphi\|_{\ld};
 \]
 therefore,
 in view of (\ref{6.2}) and the bound of type (\ref{c2b}) for $N_\delta$, we derive
\begin{equation}\label{6.3}
 \|\e^{m+1}N_\delta^\e D^{\alpha+\delta}u^{,\e}\|_{\ld}\le C\e^m\|f\|_{\ld},\quad m\ge 2,
\end{equation}
 where the constant $C$ depends on 
  $d$ and the  constants in (\ref{2}). We have just proved more than it was required, that is, 
  the
  $\ld$-bound for the term (\ref{6.1}); but we have found the majorant   of special structure for the latter (see (\ref{6.3})).
  Here is the first place we need to employ the iteration of Steklov`s smoothing operator while  
  mere Steklov`s smoothing is not enough.
  
  The differentiation of order $m$ of the  correctors in (\ref{2.7}) produces also  the 
  products
   other than those of (\ref{6.1}), namely, 
 \begin{equation}\label{6.4}
  \e^{m-|\mu|}(D^\mu N_\gamma)^\e D^{\gamma+\alpha-\mu}u^{,\e},\quad 0\le |\mu|\le m ,\,|\a|=|\gamma|=m,
  \end{equation}
  and
  \begin{equation}\label{6.5}
  \e^{m+1-|\mu|}(D^\mu N_\delta)^\e D^{\delta+\alpha-\mu}u^{,\e},\quad 0<|\mu|\le m ,\,|\a|=m,\, |\delta|=m+1.
  \end{equation}
  We handle the terms of types (\ref{6.4}) and (\ref{6.5}) by straightforward applying of Lemma \ref{LemM1}. We omit details which are similar to those used to obtain (\ref{6.3}).
 
{\bf 6.2.} We pass to the proof of (\ref{2.10}).
The  discrepancy  of the function (\ref{2.6}) in  equation (\ref{1})
admits the representation
\begin{equation}\label{6.6}
\ds{
A_\e \tilde{u}^\e{+}\tilde{u}^\e{-}f=(A_\e{+}1)\tilde{u}^\e{-}(\hat{A}_\e{+}1) u^{,\e}{+}(f^{,\e}{-}f)}
\atop\ds{
=
(A_\e \tilde{u}^\e{-}\hat{A}_\e u^{,\e}){+}
(\tilde{u}^\e{-} u^{,\e}){+}(f^{,\e}{-}f),}
\end{equation}%
where at the first step we take into account the equality
\[
(\hat{A}_\e+1) u^{,\e}=f^{,\e},\quad f^{,\e}=\Theta^\e f,
\]
obtained  by applying the smoothing operator $\Theta^\e$ to the both parts of (\ref{2.2})
 if the notation  $\Theta^\e\hat{u}^\e=u^{,\e}$ is used (see (\ref{2.8})).
 
 Comparing (\ref{2.6}) with (\ref{5.1}) and, in parallel, (\ref{6.6}) with (\ref{5.2}), we see that $\tilde{u}^\e$ correlates to  $u^{,\e}$ in (\ref{2.6}) as $v^\e$ to $u$ in (\ref{5.1}). The  calculations accomplished for the pair
$v^\e$ and $u$  in Section 5 can be repeated for the pair $\tilde{u}^\e$ and  $u^{,\e}$. By this process all expressions and all transitions make sense, which can be justified 
with the help of arguments similar to those used in Subsection 6.1. We  rely in an essential way on lemmas \ref{LemM1} and \ref{LemM2}. Thus, similarly as in (\ref{5.12}) (see also (\ref{5.2})), we obtain 
\begin{equation}\label{6.7}
\ds{
(A_\e{+}1)\tilde{u}^\e{-}(\hat{A}_\e{+}1) u^{,\e}}
\atop\ds{
=
\sum_j \e^{n_j}b_j(x/\e)\Phi_j(x)+
\sum_{|\a|=m}D^\a \sum_j \e^{n_j}\tilde{b}_j(x/\e)\tilde{\Phi}_j(x)=:F_\e,\, n_j{\ge} 2,}
\end{equation}
with the only difference that
 the functions  $\Phi_j$  and $\tilde{\Phi}_j$ in (\ref{6.7})
  coincide  with 
 the function $u^{,\e}$ (defined in (\ref{2.8})) or the derivatives $D^\nu u^{,\e}$ up to the order $2m+1$.
 The right-hand side of (\ref{6.7}) denoted by $F_\e$ has the estimate
 \begin{equation}\label{6.8}
 \|F_\e\|_{H^{-m}(\rd)}\le C\e^2\|f\|_{\ld}.
 \end{equation}
  
 As for the remaining term $f^{,\e}{-}f$  in the right-hand side of
 (\ref{6.6}), it also  possesses the desired estimate, due to (\ref{m.7}), because, for $m\ge 2$, it holds
\begin{equation}\label{6.9}
 \|f^{,\e}{-}f\|_{H^{-m}(\rd)}\le  \|f^{,\e}{-}f\|_{H^{-2}(\rd)}\le C\e^2\|f\|_{\ld}.
\end{equation}
 
 From  (\ref{6.6}) and the last two estimates, we readily come to  (\ref{2.10}). In fact,
 \[A_\e \tilde{u}^\e{+}\tilde{u}^\e{-}f=A_\e \tilde{u}^\e{+}\tilde{u}^\e-(A_\e {u}^\e{+}{u}^\e)=(A_\e+1)(\tilde{u}^\e-{u}^\e);
 \]
 therefore, using the notation from (\ref{6.7}), we have
 \[
 (A_\e+1)(\tilde{u}^\e-{u}^\e)=F_\e+(f^{,\e}{-}f),
 \]
 and it suffices to recall the energy estimate for this equation together with (\ref{6.8}) and (\ref{6.9}).

\bigskip
{\bf 6.3.} We 
conclude with some remarks.
\begin{rem}\label{remark6.1} \rm
The approximation (\ref{2.6}) can be taken with arbitrary smoothing operator (\ref{m.11}), provided that the smoothing kernel is even and, at least, piecewise $C^1$-smooth. In this case, the estimate (\ref{2.10}) remains valid.
\end{rem}
\begin{rem}\label{remark6.2}\rm
 There is another approach to prove operator type estimates in homogenization. It relies on the 
Floquet--Bloch transform and spectral expansions; thus, it is often called the spectral method. This method  is tightly linked with periodic setting because of the 
Floquet--Bloch transform which fits only this setting.
One of the first applications of the spectral approach
for  homogenization  estimates can be found in \cite{Zh89} where diffusion problems were studied. As for
high order operators, this method was applied, 
e.g., in  \cite{V},  \cite{KS},  and \cite{SS}.
\end{rem}
\begin{rem}\label{remark6.3} \rm
The method demonstrated to prove estimates   (\ref{2.10}) and  (\ref{2.13}) can be extended to   more general situations  when operators are matrix-valued. 
We choose to consider here the 
scalar problem 
only for brevity of exposition.
\end{rem}

\end{document}